\documentclass{amsart}

\newtheorem{theorem}{Theorem}[section]
\newtheorem{proposition}[theorem]{Proposition}

\newtheorem{note}[theorem]{Note}
\newtheorem{corollary}[theorem]{Corollary}

\theoremstyle{definition}
\newtheorem{definition}[theorem]{Definition}

\theoremstyle{remark}
\newtheorem{remark}[theorem]{Remark}

\numberwithin{equation}{section}

\usepackage{graphicx}
\usepackage{subfig}
\usepackage{pictexwd, dcpic}
\usepackage{amsmath, amstext}
\usepackage{amssymb}
\usepackage{amscd}
\usepackage{bm}   % for bold math
\usepackage{mathptmx}  % for postscript fonts -- requested by springer
\usepackage{hyperref}

\usepackage{color}

\definecolor{darkgreen}{rgb}{0,0.5,0}

\newcommand{\newword}[1]{\textbf{#1}}

\DeclareMathOperator{\rot}{\mathrm{rot}}
\DeclareMathOperator{\PLog}{\mathrm{PLog}}
\DeclareMathSymbol{\Rb}{\mathbin}{AMSb}{"52}
\DeclareMathSymbol{\Zb}{\mathbin}{AMSb}{"5A}

\title{Conjugacy and Dynamics in Thompson's Groups}

\author{James Belk} 
\address{Mathematics Program, Bard College, Annandale-on-Hudson, NY 12504, USA}
\email{belk@bard.edu}
\thanks{The first author gratefully acknowledges
partial support from an NSF Postdoctoral Research Fellowship while he was at Texas A\&M 
University.}

\author{Francesco Matucci}
\address{Department of Mathematics, University of Virginia, Charlottesville, VA 22904, USA}
\email{fm6w@virginia.edu}
\thanks{This work is part of the second author's Ph.D.\ thesis at Cornell University.
The second author gratefully acknowledges the Centre de Recerca Matem\`atica (CRM)
and its staff for the support received during the development of this work.}

\begin{document}

%\date{}
\maketitle

\begin{abstract}
We give a unified solution to the conjugacy problem for Thompson's groups $F$,
$T$, and $V$.  The solution uses ``strand diagrams'', which are similar in spirit to braids and generalize tree-pair diagrams for elements of Thompson's groups.
Strand diagrams are closely related to piecewise-linear functions for elements of Thompson's groups, and we use this correspondence to investigate the dynamics of elements of~$F$.
Though many of the results in this paper are known, our approach is new, and it yields elegant proofs of several old results.
\end{abstract}

\section{Introduction}

In the late 1960's, Richard J.~Thompson introduced three groups $F$, $T$, and $V$, now known collectively as \newword{Thompson's groups}.  The group $F$ consists of all piecewise-linear homeomorphisms of $[0,1]$ with finitely many breakpoints
satisfying the following conditions:
\begin{enumerate}
\item Every slope is a power of two, and\smallskip
\item Every breakpoint has dyadic rational coordinates.
\end{enumerate}
The groups $T$ and $V$ are defined similarly, except that $T$ is a group of homeomorphisms of the circle, and $V$ is a group of homeomorphisms of the Cantor set.  These three groups have been studied extensively, and we will assume that they are somewhat familiar to the reader.  For a thorough introduction, see the notes by Cannon, Floyd and Parry~\cite{cfp}.

In this paper, we present a unified solution to the conjugacy problems for Thompson's groups $F$, $T$, and~$V$.  The conjugacy problems for $F$ and $V$ have been solved separately before: a solution for $F$ was given by Guba and Sapir~\cite{gusa1} in the more general context of diagram groups, while a solution for $V$ was given by Higman~\cite{Hig} and again by Salazar-Diaz~\cite{saldiazthesis}.  As far as we know, the solution to the conjugacy problem for Thompson's group $T$ is new.  In addition, our methods work the same way for all three groups, and may be helpful for solving the conjugacy problem for certain related groups, such as the generalized Thompson groups $F(p)$ defined by Brown~\cite{brown3}, the universal central extension $\widetilde{T}$ of $T$~\cite{belkthesis}, or the braided Thompson group $BV$ defined by Brin~\cite{brin3}.

Our solution uses \newword{strand diagrams} for elements of Thompson's groups.  These are a generalization of the standard tree-pair diagrams, and were used in~\cite{belkthesis} to construct an Eilenberg-MacLane space for Thompson's group~$F$. Roughly speaking, a strand diagram is similar to a braid, but instead of twists a strand diagram has splits and merges.  In the case of~$F$, strand diagrams are closely related to the ``diagrams'' of Guba and Sapir~\cite{gusa1,gusa2}, and are essentially equivalent to the ``pictures'' of Bogley and Pride~\cite{BoPr,Pr1,Pr2}.

In addition to solving the conjugacy problem, our invariants yield information about the dynamics of elements of Thompson's groups.  This arises from an explicit correspondence between strand diagrams and piecewise-linear functions, which we describe in Section~\ref{sec:dynamics}.  Using this correspondence, we obtain a complete understanding of the dynamics of elements of $F$, giving simple proofs of several previously known results.  A more thorough investigation of dynamics using strand diagrams can be found in~\cite{matuccithesis}.

\section{Strand Diagrams \label{sec:unified-pow}}

In this section we introduce strand diagrams, which will be our main tool for solving the conjugacy problem.  We begin by realizing $F$, $T$, and $V$ as groups of strand diagrams, and then continue by defining the annular, toral, and closed strand diagrams that will serve as our conjugacy invariants.

\subsection{Strand Diagrams for $F$, $T$, and $V$ \label{sec:definition-of-strand}}

%\begin{figure}[t]
%\centering
%\begin{minipage}{1.5in}
%\centering
%\includegraphics{Fig1}
%\caption{A strand diagram.}
%\label{fig:FStrandDiagram}
%\end{minipage}
%\qquad
%\begin{minipage}{2.8in}
%\centering
%\vspace{2.4pt}\fbox{
%\includegraphics{Fig2}
%\qquad
%\includegraphics{Fig3}
%}
%\vspace{2.3pt}\caption{A split and a merge.}
%\label{fig:SplitMergePictures}
%\end{minipage}
%\end{figure}

\begin{figure}[t]
\centering
\subfloat[]{\includegraphics{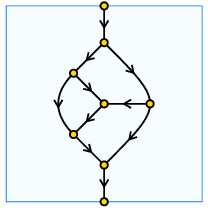}}
\qquad\qquad\qquad
\subfloat[]{
\includegraphics{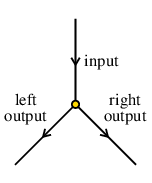}
\qquad \;
\includegraphics{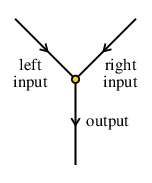}
}
\caption{(a) A strand diagram. (b) A split and a merge.}
\label{fig:FStrandDiagram}
\end{figure}

\begin{definition}A \newword{strand diagram} (see Figure~\ref{fig:FStrandDiagram}a) is a finite acyclic digraph embedded in the unit square $[0,1]\times[0,1]$, with the following properties:
\begin{enumerate}
\item The graph has a single univalent source, which lies on the top edge of the square, and a single univalent sink, which lies on the bottom edge.\smallskip
\item Every other vertex is trivalent, and is either a \newword{split} or a \newword{merge}, as shown in Figure~\ref{fig:FStrandDiagram}b.
\end{enumerate}
By convention, isotopic strand diagrams are considered equal.
\end{definition}

\begin{figure}[b]
\centering
\includegraphics{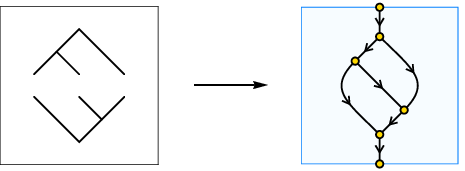}
\caption{Constructing a strand diagram from a tree-pair diagram.}
\label{fig:TreeToStrand}
\end{figure}
Strand diagrams are a generalization of tree-pair diagrams for elements of Thompson's group~$F$.  Specifically, given any tree-pair diagram, we can construct a corresponding strand diagram by gluing together the leaves of the two trees, as shown in Figure~\ref{fig:TreeToStrand}.  Not every strand diagram can be obtained in this fashion, so we can view tree-pair diagrams as a proper subset of strand diagrams.

\begin{figure}[b]
\centering
\includegraphics{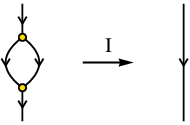}
\qquad\qquad\qquad\qquad
\includegraphics{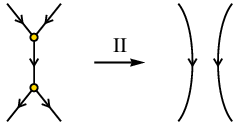}
\caption{The two reductions.}
\label{fig:Reductions}
\end{figure}
\begin{definition}A \newword{reduction} of a strand diagram is either of the two moves shown in Figure~\ref{fig:Reductions}.  A strand diagram is \newword{reduced} if it is not subject to any reductions.  Two strand diagrams are \newword{equivalent} if one can be obtained from the other by a sequence of reductions and inverse reductions.
\end{definition}

Note that two tree-pair diagrams for the same element of $F$ yield equivalent strand diagrams.  In particular, we can use type~I reductions to eliminate redundant pairs of carets that have been glued together.

\begin{proposition}Every strand diagram is equivalent to a unique reduced strand diagram.
\end{proposition}
\begin{proof}This is a straightforward application of the theory of abstract rewriting systems.\footnote{A weaker version of the Newman's Diamond Lemma suffices: see Theorem 1 in~\cite{new1}.}
In particular, the process of reduction is terminating since each reduction decreases the number of vertices, and it is straightforward to check that reductions
are locally confluent.
\end{proof}

Every reduced strand diagram can be obtained by gluing together the trees of a reduced tree-pair diagram; thus, this proposition gives us a one-to-one correspondence between equivalence classes of strand diagrams and elements of Thompson's group~$F$.
We will make this correspondence more explicit in Section~\ref{sec:dynamics} when we investigate the relationship between strand diagrams and piecewise-linear functions.

\begin{figure}[t]
\centering
\begin{tabular}{cccc}
\includegraphics{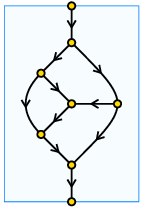} &
\includegraphics{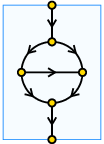} &
\includegraphics{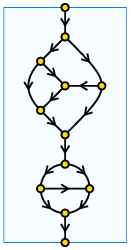} &
\includegraphics{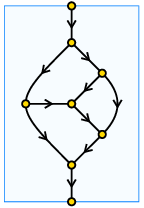} \\
$f$ & $g$ & $g\circ f$ & $g\circ f$ (reduced)
\end{tabular}
\caption{Composing elements using strand diagrams.}
\label{fig:Concatenation}
\end{figure}

The advantage of strand diagrams over tree-pair diagrams is that strand diagrams are easier to compose:

\begin{definition}The \newword{concatenation} of two strand diagrams is obtained by gluing the sink of the first to the source of the second, and then eliminating the resulting bivalent vertex (see Figure~\ref{fig:Concatenation}).
\end{definition}

Concatenation of strand diagrams corresponds to composition of elements of~$F$, as shown in Figure~\ref{fig:Concatenation}.  Note that the result is usually not reduced, so in practice we compose two elements by concatenating the strand diagrams and then reducing the result.

The following proposition should be clear from the above discussion:

\begin{proposition}The set of all reduced strand diagrams forms a group under the operation of concatenation followed by reduction. This is isomorphic to Thompson's group~$F$.\qed
\end{proposition}

\begin{figure}[b]
\centering
\includegraphics{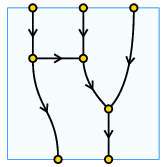}
\caption{A $(3,2)$-strand diagram.}
\label{fig:FStrandDiagram32}
\end{figure}
For technical reasons, we will also need to consider strand diagrams with more than one source or sink:

\begin{definition}Let $m$ and $n$ be positive integers.  An \newword{(\textit{m},\,\textit{n})-strand diagram} is similar to a strand diagram, except that it has $m$ univalent sources along the top of the square, and $n$~univalent sinks along the bottom of the square (see Figure~\ref{fig:FStrandDiagram32}).
\end{definition}

Note that a ``strand diagram'' is the same thing as a $(1,1)$-strand diagram.  As with strand diagrams, every $(m,n)$-strand diagram is equivalent to a unique reduced $(m,n)$-strand diagram.  Moreover, every reduced $(m,n)$-strand diagram can be obtained by attaching together the leaves of two binary forests, one with $m$~trees and the other with $n$~trees.

The set of all reduced $(m,n)$-strand diagrams forms a groupoid over the positive integers.  In particular, we can concatenate any $(i,j)$-strand diagram with any $(j,k)$-strand diagram by attaching the sinks of the first to the sources of the second.  The isotropy group of this groupoid at the point~$1$ (or indeed at any point) is isomorphic to Thompson's group~$F$.

We now turn to strand diagrams for Thompson's groups $T$ and $V$.  For Thompson's group~$T$, we can use the following class of diagrams:

\begin{figure}[t]
\centering
\subfloat[]{\includegraphics{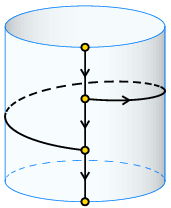}}
\qquad\qquad\qquad
\subfloat[]{\includegraphics{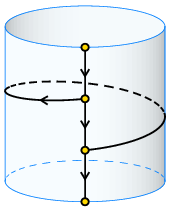}}
\qquad\qquad\qquad
\subfloat[]{\includegraphics{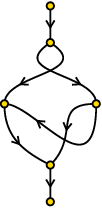}}
\caption{(a) A cylindrical strand diagram. (b) An isotopic cylindrical strand diagram, obtained from the first by a Dehn twist. (c) An abstract strand diagram.}
\label{fig:CylStrandDiagram}
\end{figure}
\begin{definition} A \newword{cylindrical strand diagram} (see Figure~\ref{fig:CylStrandDiagram}a) is a finite acyclic digraph embedded in the cylinder $S^1 \times [0,1]$, with the following properties:
\begin{enumerate}
\item The graph has a single univalent source, which lies on the top circle of the cylinder, and a single univalent sink, which lies on the bottom circle of the cylinder.\smallskip
\item Every other vertex is either a split or a merge.
\end{enumerate}
Isotopic cylindrical strand diagrams are considered equal.
\end{definition}

Note that the source and sink of a cylindrical strand diagram may rotate around their respective circles during an isotopy.  Therefore, two cylindrical strand diagrams that differ by a Dehn twist of the cylinder are considered equal (see Figure~\ref{fig:CylStrandDiagram}b).

Note also that we can glue two binary trees together on the cylinder via any cyclic permutation of the leaves.
This lets us view tree-pair diagrams for elements of~$T$ as a proper subset of cylindrical strand diagrams.

It is less obvious how to define strand diagrams for Thompson's group~$V$.  For the following definition, recall that a \newword{rotation system} on a graph is an assignment of a circular order to the edges incident on each vertex~\cite{moharthomassen}.

\begin{definition} An \newword{abstract strand diagram} is a finite acyclic digraph together with a rotation system, having the following properties:
\begin{enumerate}
\item The graph has a single univalent source and a single univalent sink.\smallskip
\item Every other vertex is either a split or a merge.
\end{enumerate}
Two abstract strand diagrams are considered equal if there exists a digraph isomorphism between them that is compatible with the corresponding rotation systems.
\end{definition}

Figure~\ref{fig:CylStrandDiagram}c shows an abstract strand diagram drawn in the plane.  By convention, the rotation system is indicated by the counterclockwise order of the edges at each vertex.  Note that the rotation system on an abstract strand diagram lets us define the left and right outputs for a split, and the left and right inputs for a merge.

Given a tree-pair diagram for an element of~$V$, we can glue the leaves of the trees together along the indicated bijection to obtain an abstract strand diagram.  This lets us view tree pair diagrams for elements of~$V$ as a subset of abstract strand diagrams.

All of the constructions we defined for strand diagrams can be extended to cylindrical and abstract strand diagrams.  In particular:
\begin{figure}[t]
\centering
\includegraphics{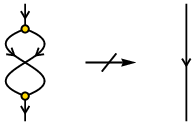}
\qquad\qquad\qquad\qquad
\includegraphics{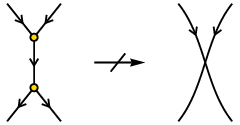}
\caption{Invalid reductions.}
\label{fig:Nonreductions}
\end{figure}
\begin{enumerate}
\item We can reduce either type of diagram using the two reduction moves shown in Figure~\ref{fig:Reductions}.  When applying these moves, it is very important to respect the counterclockwise orientation of the edges shown in Figure~\ref{fig:Reductions}.  In particular, neither of the moves shown in Figure~\ref{fig:Nonreductions} is a valid reduction.

    \rule{0pt}{0pt}Assuming we apply only valid reduction moves, every cylindrical strand diagram is equivalent to a unique reduced cylindrical strand diagram, and every abstract strand diagram is equivalent to a unique reduced abstract strand diagram.\smallskip

\item We can concatenate two cylindrical strand diagrams or two abstract strand diagrams.  Under the operation of concatenation followed by reduction, reduced cylindrical strand diagrams form a group isomorphic to Thompson's group~$T$, and reduced abstract strand diagrams form a group isomorphic to Thompson's group~$V$. \smallskip

\begin{figure}[b]
\centering
\subfloat[]{\includegraphics{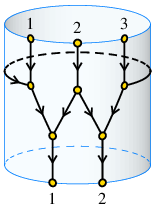}}
\qquad\qquad\qquad\qquad
\subfloat[]{\includegraphics{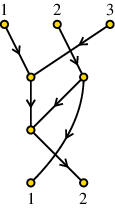}}
\caption{(a) A cylindrical $(3,2)$-strand diagram. (b) An abstract $(3,2)$-strand diagram.}
\label{fig:CylAbs32Diagrams}
\end{figure}
\item We can also define groupoids of cylindrical $(m,n)$-strand diagrams and abstract $(m,n)$-strand diagrams, though there is a slight complication.  If we wish to be able to concatenate $(m,n)$-diagrams unambiguously, we must include numberings of the sources and sinks, as shown in Figure~\ref{fig:CylAbs32Diagrams}.  For a cylindrical diagram, these numbers are required to appear in counterclockwise order around each circle.  For an abstract diagram, the sources and sinks can be numbered in any order.  When concatenating two diagrams, we use these numbers to determine which sources should be attached to which sinks.
\end{enumerate}

\subsection{Conjugacy Invariants}
Given a strand diagram in the unit square, we can \newword{close}
it to obtain a graph embedded in an annulus, as shown in Figure~\ref{fig:ClosingStrandDiagram}. More generally, we can close any $(m,n)$-strand diagram for which~$m=n$.  This prompts the following definition:
\begin{figure}[t]
\centering
\includegraphics{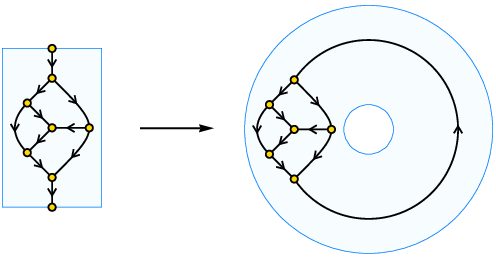}
\caption{Closing a strand diagram to obtain an annular strand diagram.}
\label{fig:ClosingStrandDiagram}
\end{figure}

\begin{definition}An \newword{annular strand diagram} is a finite directed topological graph embedded in the annulus $[0,1]\times S^1$, with the following properties:
\begin{enumerate}
\item Every vertex is either a split or a merge.\smallskip
\item Every directed cycle winds counterclockwise around the central hole.
\end{enumerate}
As with strand diagrams, isotopic annular strand diagrams are considered equal.
\end{definition}

In this definition, a \newword{topological graph} refers to a graph without bivalent vertices, which may also contain one or more \newword{free loops} (closed loops without any vertices).  Note that a free loop counts as a directed cycle, so by property~(2) every free loop in an annular strand diagram must wind counterclockwise around the central hole.

Though the underlying spaces are homeomorphic, an annular strand diagram is actually very different from a cylindrical strand diagram.  In a cylindrical strand diagram, the horizontal direction is a circle, but there are still sources and sinks along the top and bottom. In an annular strand diagram, it is the vertical direction that has been made into a circle. We shall consistently distinguish between these cases by the order of the factors in a product: $S^1\times[0,1]$ for the cylinder, and $[0,1]\times S^1$ for the annulus.

\begin{figure}[b]
\centering
\fbox{\includegraphics{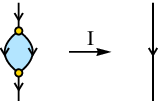}}
\hfill
\fbox{\includegraphics{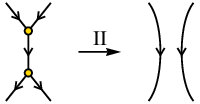}}
\hfill
\fbox{\includegraphics{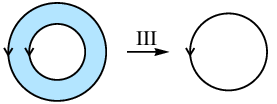}}
\caption{The three reduction moves for an annular strand diagram.  In the first move, the shaded region must be a topological disk.  In the third move, both loops must be free loops, and the region between must be a topological annulus that does not contain any vertices.}
\label{fig:ClosedReductions}
\end{figure}

\begin{figure}[t]
\centering
\includegraphics{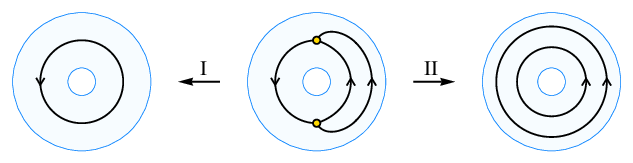}
\caption{The annular strand diagram in the center can be reduced using either a type~I move or a type~II move.  A type~III move is required to reconcile the results.}
\label{fig:Type3ReductionNeeded}
\end{figure}
We can reduce annular strand diagrams using the three reduction moves shown in Figure~\ref{fig:ClosedReductions}.  The third move is necessary to make reductions of annular strand diagrams locally confluent.  In particular, Figure~\ref{fig:Type3ReductionNeeded}
shows an annular strand diagram which can be reduced using either a type~I move or a type~II move, but for which the two results can only be reconciled using a type~III move.  As long as we use all three reductions, though, every annular strand diagram is equivalent to a unique reduced annular strand diagram.

The following theorem explains our interest in annular strand diagrams.  We will prove this theorem in Section~\ref{sec:characterize-conjugacy}.

\begin{theorem}Let $f$ and $g$ be elements of Thompson's group $F$.  Let\/ $\mathfrak{f}$ and\/ $\mathfrak{g}$ be strand diagrams for $f$ and $g$, and let\/ $\mathfrak{f}'$ and\/ $\mathfrak{g}'$ be the reduced annular strand diagrams obtained by closing\/ $\mathfrak{f}$ and\/ $\mathfrak{g}$ and reducing.  Then $f$ and $g$ are conjugate if and only if\/ $\mathfrak{f}'$ and\/ $\mathfrak{g}'$ are isotopic.
\label{thm:ConjugacyF}
\end{theorem}

This theorem provides a geometric solution
to the conjugacy problem in Thompson's group~$F$.  (Recall that isotopy of graph embeddings can be checked algorithmically using dual graphs.) This procedure is surprisingly fast: it is shown in \cite{hos} that an algorithm based on the above theorem can be implemented in linear time.

Theorem~\ref{thm:ConjugacyF} is very similar to the solution to the conjugacy
problem given by Guba and Sapir in~\cite{gusa1}, using ``diagrams'' instead of strand
diagrams.  Indeed, our strand diagrams
are essentially just the dual graphs to the
diagrams that they use.
%\red{Removed footnote. We talk about the Bogley-Pride citation in the intro only}
%\footnote{
%\blue{Strand diagrams are also essentially equivalent to the ``pictures'' used in the work of Bogley and Pride.  See \cite{BoPr},  \cite{Pr1}, and~\cite{Pr2}.}}
The advantage of strand diagrams is that they generalize
easily to Thompson's groups $T$ and~$V$.

\begin{figure}[b]
\centering
\subfloat[]{\includegraphics{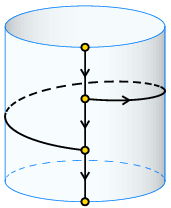}}
\hfill
\subfloat[]{\includegraphics{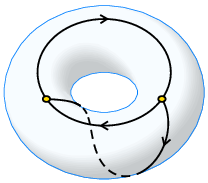}}
\hfill
\subfloat[]{\includegraphics{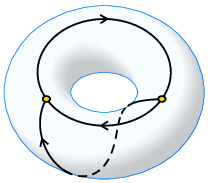}}
\caption{(a) A cylindrical strand diagram. (b) The corresponding toral strand diagram. (c)~The same toral strand diagram after a Dehn twist around $S^1\times\{1\}$.}
\label{fig:ToralStrandDiagram}
\end{figure}
We begin by describing the generalization to Thompson's group~$T$.  Given a cylindrical strand diagram, we can \newword{close} it by gluing together the two circles of the cylinder, identifying the source and the sink (see Figure~\ref{fig:ToralStrandDiagram}).  The result is a graph embedded on a torus $S^1 \times S^1$.

For the following definition, let $c$ denote the cohomology class in $H^1(S^1\times S^1)$ that measures the intersection number with the circle $S^1\times\{1\}$ (the circle obtained by gluing together the top and bottom circles of the cylinder).  Note that this cohomology class is invariant under the Dehn twist around the circle $S^1\times\{1\}$, and indeed the stabilizer of $c$ in the orientation-preserving mapping class group is precisely the cyclic subgroup generated by this Dehn twist.

\begin{definition}A \newword{toral strand diagram} is a finite directed topological graph embedded on the torus $S^1\times S^1$, with the following properties:
\begin{enumerate}
\item Every vertex is either a split or a merge.\smallskip
\item Every directed cycle evaluates to a positive value under~$c$.
\end{enumerate}
Two toral strand diagrams are considered equal if they are isotopic, or if they differ by finitely many Dehn twists around the circle $S^1\times\{1\}$ (see Figure~\ref{fig:ToralStrandDiagram}c).
\end{definition}

We can reduce toral strand diagrams in exactly the same way that we reduce annular
strand diagrams, and each toral strand diagram is equivalent to a unique reduced toral
strand diagram.  The following theorem provides a solution to the conjugacy problem
in Thompson's group~$T$.  We will prove this theorem in
Section~\ref{sec:characterize-conjugacy}.

\begin{theorem}Let $f$ and $g$ be elements of Thompson's group $T$.  Let\/ $\mathfrak{f}$ and\/ $\mathfrak{g}$ be cylindrical strand diagrams for $f$ and $g$, and let\/ $\mathfrak{f}'$ and\/ $\mathfrak{g}'$ be the reduced toral strand diagrams obtained by closing\/ $\mathfrak{f}$ and\/ $\mathfrak{g}$ and reducing.  Then $f$ and $g$ are conjugate if and only if\/ $\mathfrak{f}'$ and\/ $\mathfrak{g}'$ are equal.
\label{thm:ConjugacyT}
\end{theorem}

Note that checking whether $\mathfrak{f}'$ and $\mathfrak{g}'$ are equal involves checking whether they differ by finitely many Dehn twists around $S^1\times\{1\}$.  This can be determined algorithmically as follows:
\begin{enumerate}
\item Enumerate all digraph isomorphisms between the two diagrams.\smallskip
\item For each digraph isomorphism, use dual graphs to check whether the isomorphism extends to an orientation-preserving self-homeomorphism of the torus.\smallskip
\item For each isomorphism that extends, check whether the corresponding homeomorphism preserves the cohomology class~$c$.  If it does, then the two toral strand diagrams are equal.  If no such isomorphism is found, then the two toral strand diagrams are different.
\end{enumerate}
Here $c$ is the conjugacy class that measures winding number around the central hole.
%The algorithm above is only one possibility for
%checking equality---see Section~??? \red{[We need to clear
% this up.]} for an alternative algorithm for testing the equality of
% two reduced toral strand diagrams.

Finally, we would like to discuss the solution to the conjugacy problem in Thompson's group~$V$.  Because abstract strand diagrams for $V$ do not lie on a surface, we must replace the ``central hole'' of the annulus or torus with a cohomology class:

\begin{definition}An \newword{abstract closed strand diagram} is an ordered pair $(\Gamma,c)$, where $\Gamma$ is a finite directed topological graph with a rotation system, and $c\in H^1(\Gamma;\mathbb{Z})$, having the following properties:
\begin{enumerate}
\item Every vertex of $\Gamma$ is either a split or a merge.\smallskip
\item For every directed cycle $\alpha$ in $\Gamma$, we have $c(\alpha) > 0$.
\end{enumerate}
Two abstract closed strand diagrams $(\Gamma,c)$ and $(\Gamma',c')$ are considered equal if there exists an isomorphism $\varphi\colon\Gamma\to\Gamma'$ of directed graphs with rotation systems such that $\varphi^*(c') = c$.
\end{definition}

\begin{figure}[t]
\centering
\includegraphics{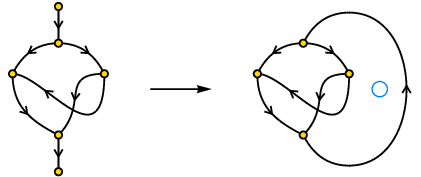}
\caption{Closing an abstract strand diagram.}
\label{fig:AbstractClosure}
\end{figure}
We can \newword{close} any abstract strand diagram by attaching the source and the sink, resulting in an abstract closed strand diagram (see Figure~\ref{fig:AbstractClosure}).  The cohomology class $c$ is defined by the $1$-cochain that evaluates to $1$ on the new edge obtained by gluing the source to the sink, and $0$ on every other edge.

\begin{figure}[b]
\centering
\includegraphics{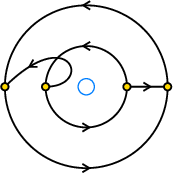}
\qquad\qquad\qquad\qquad
\includegraphics{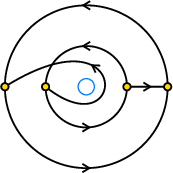}
\caption{Two abstract closed strand diagrams with slightly different cutting classes.}
\label{fig:AbstractClosedStrandDiagrams}
\end{figure}
Figures~\ref{fig:AbstractClosure} and \ref{fig:AbstractClosedStrandDiagrams} show drawings of several abstract closed strand diagrams.  In each drawing, the rotation system~$r$ is indicated by the counterclockwise order of the edges around each vertex, and the cohomology class $c$ corresponds to counterclockwise winding number around the indicated ``hole''.  It is always possible to draw an abstract closed strand diagram in this way, since the elements of $H^1(\Gamma;\mathbb{Z})$ are in one-to-one correspondence with homotopy classes of maps from $\Gamma$ to the punctured plane.

We shall refer to the cohomology class $c$ as the \newword{cutting class} of the abstract closed strand diagram.  Note that the cutting class is an essential part of the invariant, and is not determined by the underlying graph $\Gamma$.  For example, Figure~\ref{fig:AbstractClosedStrandDiagrams} shows two abstract closed strand diagrams with isomorphic graphs and rotation systems, but slightly different cutting classes.  These two diagrams are not equal, and correspond to different conjugacy classes in Thompson's group~$V$.

\begin{remark}
Though we can reduce an abstract closed strand diagram using the three reductions shown in Figure~\ref{fig:ClosedReductions}, some care must be taken in interpreting these moves.  In particular, the following restrictions apply:
\begin{enumerate}
\item For a type~I reduction, the closed cycle $\alpha$ formed by the two parallel edges must satisfy~$c(\alpha)=0$. This replaces the requirement that the two parallel edges bound a disk.
\smallskip
\item For a type~III reduction, the two free loops $\beta$ and $\gamma$ must satisfy $c(\beta)=c(\gamma)$. This replaces the requirement that the two free loops bound an annulus.
\end{enumerate}
In addition, the counterclockwise order of the edges at each vertex must be exactly as shown in Figure~\ref{fig:ClosedReductions}.  That is, the two moves shown in Figure~\ref{fig:Nonreductions} are not allowed.

Note also that, for each type of reduction, the cutting class on the original diagram induces a cutting class on the reduced diagram in a natural way.
\label{rmk:ReduceClosedAbstract}
\end{remark}

The following theorem characterizes conjugacy in Thompson's group~$V$.  We will prove it in Section~\ref{sec:characterize-conjugacy}.

\begin{theorem}Let $f$ and $g$ be elements of Thompson's group $V$.  Let\/ $\mathfrak{f}$ and\/ $\mathfrak{g}$ be abstract strand diagrams for $f$ and $g$, and let\/ $\mathfrak{f}'$ and\/ $\mathfrak{g}'$ be the reduced abstract closed strand diagrams obtained by closing\/ $\mathfrak{f}$ and\/ $\mathfrak{g}$ and reducing.  Then $f$ and $g$ are conjugate if and only if\/ $\mathfrak{f}'$ and\/ $\mathfrak{g}'$ are equal.
\label{thm:ConjugacyV}
\end{theorem}

Note that we can check equality of abstract closed strand diagrams algorithmically by enumerating all isomorphisms between the directed graphs and then checking whether the isomorphisms are compatible with the rotation systems and cutting classes.  Thus the above theorem provides a solution to the conjugacy problem in $V$.

\newcommand{\eqf}{[\,\mathfrak{f}\,]}
\newcommand{\eqg}{[\hspace{0.7pt}\mathfrak{g}\hspace{0.7pt}]}
\newcommand{\eqh}{[\hspace{0.7pt}\mathfrak{h}_1\hspace{0.5pt}]}
\newcommand{\eqhh}{[\hspace{0.7pt}\mathfrak{h}_2\hspace{0.5pt}]}

\section{Solving the Conjugacy Problem \label{sec:characterize-conjugacy}}

The goal of this section is to prove Theorems~\ref{thm:ConjugacyF}, \ref{thm:ConjugacyT}, and \ref{thm:ConjugacyV}, which solve the conjugacy problems for $F$, $T$, and~$V$.  Our approach is to handle all three cases simultaneously, pointing out the differences when necessary.

Throughout this section, we will use the term \newword{(\textit{m},\,\textit{n})-diagram} to refer to either an \mbox{$(m,n)$-strand diagram}, a cylindrical $(m,n)$-strand diagram, or an abstract $(m,n)$-strand diagram.  Similarly, we will use the term \newword{closed diagram} to refer to either an annular strand diagram, a toral strand diagram, or an abstract closed strand diagram.

Recall that two diagrams are \newword{equivalent} if they differ by a sequence of reductions and inverse reductions.  For convenience we will be working primarily in the groupoid of equivalence classes of $(m,n)$-diagrams.  If $\mathfrak{f}$ is an $(m,n)$-diagram, we will let $\eqf$ denote the corresponding element of the groupoid. Here is our main theorem:

\begin{theorem}Let\/ $\mathfrak{f}$ be an $(m,m)$-diagram, and let\/ $\mathfrak{g}$ be an $(n,n)$-diagram. Then\/ $\eqf$ and\/ $\eqg$ are conjugate if and only if the closures of\/ $\mathfrak{f}$ and\/ $\mathfrak{g}$ are equivalent.
\label{thm:MainConjugacyTheorem}
\end{theorem}

Note that Theorems~\ref{thm:ConjugacyF}, \ref{thm:ConjugacyT}, and \ref{thm:ConjugacyV} follow immediately.  In particular, two elements of $F$, $T$, and $V$ are conjugate within the group if and only if they are conjugate within the corresponding groupoid.

The proof of Theorem~\ref{thm:MainConjugacyTheorem} occupies the remainder of this section.  Most of the content of the proof is covered in the following pair of propositions.

\begin{proposition}Let\/ $\mathfrak{f}$ be an $(m,m)$-diagram, let\/ $\mathfrak{f}'$ be its closure, and let\/ $\mathfrak{g}'$ be a closed diagram obtained by applying a reduction to\/ $\mathfrak{f}'$.  Then there exists an $(n,n)$-diagram\/ $\mathfrak{g}$ whose closure is\/ $\mathfrak{g}'$ such that\/ $\eqg$ is conjugate to\/~$\eqf$.
\label{prop:ReductionConjugate}
\end{proposition}
\begin{proof}Throughout this proof, we will use the term ``gluing points'' to refer to the points in~$\mathfrak{f}'$ obtained by gluing together the sources and sinks of~$\mathfrak{f}$. There are three cases, based on the type of reduction.

\begin{figure}[b]
\centering
\subfloat[]{\includegraphics{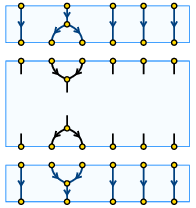}}
\hfill
\subfloat[]{\includegraphics{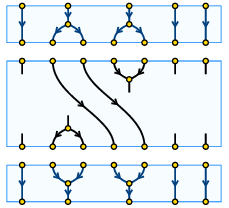}}
\hfill
\subfloat[]{\includegraphics{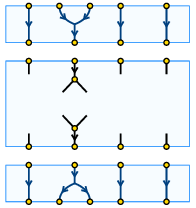}}
\caption{(a) Conjugating $\mathfrak{f}$ to perform a type~I reduction on $\mathfrak{f}'$. (b) A type~I reduction involving more than one pair of gluing points. (c) A type~II reduction.}
\label{fig:ReductionConjugate12}
\end{figure}
Suppose the reduction of $\mathfrak{f}'$ is of type~I.  If this corresponds to
a type~I reduction on $\mathfrak{f}$, then we are done.  However, it is possible that the reduction has a gluing point on each of the parallel edges, as shown in Figure~\ref{fig:ReductionConjugate12}a.  In this case, conjugating by an element that merges the two sinks will produce~$\mathfrak{g}$, as shown in the figure.  Finally, it is possible that each of the parallel edges contains several gluing points, as shown in Figure~\ref{fig:ReductionConjugate12}b.  In this case, we must conjugate by an element that merges each pair of relevant sinks. Note that the number of gluing points must be the same on the two parallel edges because of the requirement that they bound a disk (see Figure~\ref{fig:ClosedReductions})---or, in the case of abstract diagrams, because of the cohomology requirement stated in Remark~\ref{rmk:ReduceClosedAbstract}.

\begin{figure}[t]
\centering
\subfloat[]{\includegraphics{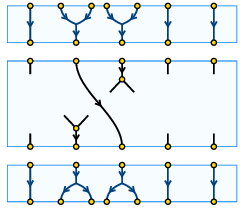}}
\hfill
\subfloat[]{\includegraphics{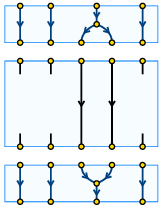}}
\hfill
\subfloat[]{\includegraphics{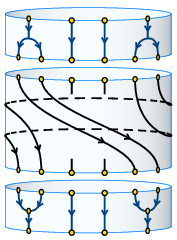}}
\caption{(a) A type~II reduction involving more than one gluing point. (b) A type~III reduction. (c) A type~III reduction involving more than one pair of gluing points.}
\label{fig:ReductionConjugate23}
\end{figure}
The argument for a type~II reduction is fairly similar.  If the reduction corresponds to a reduction on~$\mathfrak{f}$, we are done.  Otherwise we must conjugate by an element with one or more splits, as shown in Figures~\ref{fig:ReductionConjugate12}c and Figures~\ref{fig:ReductionConjugate23}a.

Finally, suppose the reduction of $\mathfrak{f}'$ is of type~III. If there is one gluing point on each free loop, then we can conjugate by an element with a single split, as shown in Figure~\ref{fig:ReductionConjugate23}b.  This is the only possibility in the annular case, but in the toral and abstract cases we might have multiple gluing points on each free loop, as shown in Figure~\ref{fig:ReductionConjugate23}c.  In this case, we must conjugate by an element with multiple splits, as shown. Again, note that the number of gluing points must be the same on each loop because of the requirement that they bound an annulus (see Figure~\ref{fig:ClosedReductions})---or, in the case of abstract diagrams, because of the cohomology requirement stated in Remark~\ref{rmk:ReduceClosedAbstract}.
\end{proof}

\begin{proposition}Let\/ $\mathfrak{f}$ be an $(m,m)$-diagram, let\/ $\mathfrak{g}$ be an $(n,n)$-diagram. If the closures of\/ $\mathfrak{f}$ and\/ $\mathfrak{g}$ are equal, then\/ $\eqf$ and\/ $\eqg$ are conjugate.
\label{prop:EqualConjuate}
\end{proposition}
\begin{proof}Let $\{\mathfrak{f}_k\}_{k\in\mathbb{Z}}$ be a collection of disjoint copies of $\mathfrak{f}$, and consider the ``infinite strand diagram'' $\mathfrak{f}^\infty$ obtained by gluing the sinks of each $\mathfrak{f}_k$ to the sources of $\mathfrak{f}_{k+1}$.  If $\mathfrak{f}$ is annular or toral, then $\mathfrak{f}^\infty$ is an infinite graph embedded in $[0,1]\times\mathbb{R}$ or $S^1\times\mathbb{R}$, respectively.  If $\mathfrak{f}$ is abstract, then $\mathfrak{f}^\infty$ is simply an infinite directed graph with a rotation system.

Let $\mathfrak{f}'$ be the closure of $\mathfrak{f}$, and note that $\mathfrak{f}^\infty$ is an infinite-sheeted cover of $\mathfrak{f'}$.  In the annular or toral case, this cover is induced by the covering map from $[0,1]\times\mathbb{R}$ to the annulus, or from $S^1\times\mathbb{R}$ to the torus.  In the abstract case, the kernel of the cutting class of $\mathfrak{f}'$ is an infinite-index subgroup of $\pi_1(\mathfrak{f}')$, and $\mathfrak{f}^\infty$ is the cover corresponding to this subgroup.  Either way, the cover is entirely determined by~$\mathfrak{f}'$.

Since the closure $\mathfrak{g}'$ of $\mathfrak{g}$ is equal to~$\mathfrak{f}'$, it follows that the infinite strand diagram $\mathfrak{g}^\infty$ must be the same as $\mathfrak{f}^\infty$.  To be precise:
\begin{enumerate}
\item If $\mathfrak{f}'$ and $\mathfrak{g}'$ are annular, then they are isotopic, and this isotopy lifts to an isotopy between $\mathfrak{f}^\infty$ and $\mathfrak{g}^\infty$.\smallskip
\item If $\mathfrak{f}'$ and $\mathfrak{g}'$ are toral, then either they are isotopic, or they differ by finitely many Dehn twists around $S^1\times\{1\}$.  An isotopy on the torus lifts to an isotopy on the infinite cylinder, and a Dehn twist around $S^1\times\{1\}$ also lifts to an isotopy on the infinite cylinder.  Either way, $\mathfrak{f}^\infty$ and $\mathfrak{g}^\infty$ are isotopic.\smallskip
\item If $\mathfrak{f}'$ and $\mathfrak{g}'$ are abstract, then there is an isomorphism between them that preserves the cutting class, and this lifts to an isomorphism between $\mathfrak{f}^\infty$ and $\mathfrak{g}^\infty$.
\end{enumerate}
Therefore, $\mathfrak{f}^\infty$ can be expressed either as the union $\bigcup_{k\in\mathbb{Z}}\mathfrak{f}_k$ of infinitely many copies of~$\mathfrak{f}$, or as the union $\bigcup_{k\in\mathbb{Z}}\mathfrak{g}_k$ of infinitely many copies of~$\mathfrak{g}$.  Moreover, there exists a deck transformation $t$ of $\mathfrak{f}^\infty$ such that $t(\mathfrak{f}_k) = \mathfrak{f}_{k+1}$ and $t(\mathfrak{g}_k) = \mathfrak{g}_{k+1}$ for each~$k$.

Now, since $\mathfrak{f}$ and $\mathfrak{g}$ are compact, there exists a sufficiently large $N\in\mathbb{N}$ so that $\mathfrak{g}_N$ is entirely contained in $\bigcup_{k>0} \mathfrak{f}_k$.  Let
\[
\mathfrak{h}_0 \;=\; \left(\,\bigcup_{k\geq 0} \mathfrak{f}_k\right) \,\cap\, \left(\,\bigcup_{k< N} \mathfrak{g}_k\right) \qquad\text{and}\qquad \mathfrak{h}_1 \;=\; \left(\,\bigcup_{k > 0} \mathfrak{f}_k\right) \,\cap\, \left(\,\bigcup_{k \leq N} \mathfrak{g}_k\right).
\]
Then $\mathfrak{h}_0$ and $\mathfrak{h}_1$
are strand diagrams
%of the same type
%as $\mathfrak{f}$ and $\mathfrak{g}$,
%\red{Same type? We had removed occurrences of same type. Shall we remove
%``of the same type as $f$ and $g$''?}
and $\mathfrak{h}_1 = t(\mathfrak{h}_0)$, so $\mathfrak{h}_0$ and $\mathfrak{h}_1$ are equal to the same strand diagram~$\mathfrak{h}$.  But $\mathfrak{f}_0 \cup \mathfrak{h}_1 = \mathfrak{h}_0 \cup \mathfrak{g}_N$, so the concatenation of $\mathfrak{f}$ with~$\mathfrak{h}$ is equal to the concatenation of $\mathfrak{h}$ with~$\mathfrak{g}$.
\end{proof}

We are now ready to prove the main theorem.

\begin{proof}[Theorem~\ref{thm:MainConjugacyTheorem}] Let $\mathfrak{f}$ be an $(m,m)$-diagram, and let $\mathfrak{g}$ be an $(n,n)$-diagram of the same type. Suppose first that $\eqf$ and $\eqg$ are conjugate.  Then $\mathfrak{f}$ is equivalent to a concatenation of the form $\mathfrak{h}^{-1}\mathfrak{g}\mathfrak{h}$ for some $(m,n)$-diagram $\mathfrak{h}$.  Then the closure of $\mathfrak{f}$ is equivalent to the closure of $\mathfrak{h}^{-1}\mathfrak{g}\mathfrak{h}$.  This is the same as the closure of $\mathfrak{g}\mathfrak{h}\mathfrak{h}^{-1}$, which is equivalent to the closure of~$\mathfrak{g}$.

For the converse, suppose that closures of $\mathfrak{f}$ and $\mathfrak{g}$ are equivalent.  Then applying reductions to either $\mathfrak{f}$ or $\mathfrak{g}$ eventually results in the same reduced closed diagram~$\mathfrak{h'}$.  By Proposition~\ref{prop:ReductionConjugate}, there exist diagrams $\mathfrak{h}_1$ and $\mathfrak{h}_2$ whose closures are $\mathfrak{h}'$ such that $\eqf$ is conjugate to $\eqh$, and $\eqg$ is conjugate to~$\eqhh$.  By Proposition~\ref{prop:EqualConjuate}, the elements $\eqh$ and $\eqhh$ must themselves be conjugate, and therefore $\eqf$ is conjugate to~$\eqg$.
\end{proof} 
\section{Structure of Reduced Closed Strand Diagrams}
\label{ssec:structure-closed-diagrams}

\begin{figure}[b]
 \centering
  \includegraphics{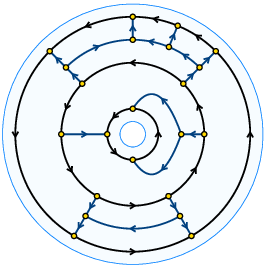}\qquad\includegraphics{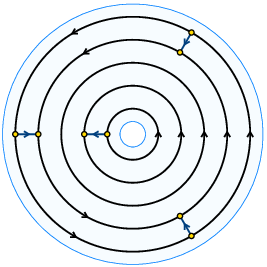}
  \caption{A pair of reduced annular strand diagrams.  The first is connected, while the second has three connected components.  The directed cycles are shown in black.}
  \label{fig:reduced-closed-example}
\end{figure}

In this section we briefly investigate the structure of reduced closed strand diagrams.  While reading the subsequent proof, it might help to refer to Figure~\ref{fig:reduced-closed-example}, which gives two examples of reduced annular strand diagrams.

As in the previous section, we will used the phrase \newword{closed diagram} to refer to either an annular strand diagram, a toral strand diagram, or an abstract closed strand diagram.  A closed strand diagram is \textbf{reduced} if it is not subject to any of the three reduction moves shown in Figure~\ref{fig:ClosedReductions}.

We will use the following terminology for certain kinds of directed cycles:
\begin{enumerate}
\item A \textbf{free loop} is a directed cycle with no vertices.\smallskip
\item A \textbf{split loop} is a directed cycle with splits, but no merges.\smallskip
\item A \textbf{merge loop} is a directed cycle with merges, but no splits.
\end{enumerate}
For example, the first annular strand diagram in Figure~\ref{fig:reduced-closed-example} has one split loop and two merge loops, while the second annular strand diagram in the figure has two split loops, two merge loops, and one free loop.

\begin{proposition}
Let\/ $\mathfrak{f}$ be any reduced closed strand diagram.  Then:
\begin{enumerate}
\item Every component of\/ $\mathfrak{f}$ has at least one directed cycle.\smallskip
\item Every directed cycle in\/ $\mathfrak{f}$ is either a free loop, a split loop, or a merge loop.\smallskip
\item Any two directed cycles in\/ $\mathfrak{f}$ are disjoint, and no directed cycle intersects itself.
\end{enumerate}
\end{proposition}
\begin{proof}For statement (1), observe that every vertex in $\mathfrak{f}$ has at least one outgoing edge.  Therefore, there exists an infinite directed path beginning at each vertex of~$\mathfrak{f}$.  Such a directed path must eventually intersect itself, and therefore $\mathfrak{f}$ must contain a directed cycle.

For statement (2), suppose to the contrary that some directed cycle of $\mathfrak{f}$ has both merges and splits.  Then at least one edge of the cycle must begin at a merge and end at a split.  Then this edge is subject to a type~II reduction, which contradicts the assumption that $\mathfrak{f}$ is reduced.

For statement~(3), observe that any two intersecting directed cycles would have to merge together and then split apart, again contradicting the assumption that $\mathfrak{f}$ is reduced.  The same reasoning shows that a directed cycle cannot intersect itself.
\end{proof}

Thus every reduced closed strand diagram consists of finitely many disjoint directed cycles, joined together by acyclic directed subgraphs.  A bit more can be said for reduced annular and toral strand diagrams:
\begin{enumerate}
\item For a reduced annular strand diagram, every component contains a directed cycle, and therefore every component must surround the central hole.  It follows that every annular strand diagram consists of finitely many concentric components, each of which is itself an annular strand diagram.  Note that a single component may have arbitrarily many directed cycles, though in a reduced diagram these cycles must alternate concentrically between merge loops and split loops.\smallskip
\item For a reduced toral strand diagram, we know that every directed cycle must have positive winding number around the central hole.  Since the directed cycles cannot intersect, all the directed cycles in a reduced toral strand diagram must be homotopic as loops in the torus.  Specifically, each directed cycle must rotate around the central hole $n$ times, and must rotate around the torus $k$ times in the other direction, where $k$ and $n$ are relatively prime.  By performing Dehn twists, we may assume that $0\leq k < n$.  This gives us a rational invariant $k/n\in \mathbb{Q}\cap [0,1)$ of the conjugacy class, which we refer to as the \textbf{rotation number}.  Though we shall not prove it here, the rotation number of a conjugacy class for $T$ is the same as the dynamical rotation number of any element of the class, viewed as a homeomorphism of the circle.
\end{enumerate}

\section{Dynamics in Thompson's group $F$ \label{sec:dynamics}}

In this section we show how the structure of an annular strand diagram is related to the dynamics of elements of the corresponding conjugacy class. For simplicity,
we restrict to the case of Thompson's group~$F$, although similar results hold for $T$ and~$V$.

\subsection{Strand diagrams as piecewise-linear functions \label{sec:representation-elements}}

\begin{figure}[t]
\centering
\subfloat[]{\includegraphics{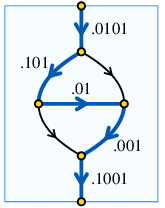}}
\qquad\qquad
\subfloat[]{\includegraphics{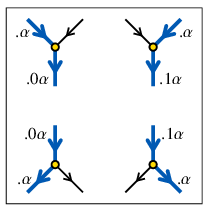}}
\caption{(a) Computing $f(.0101)$ using a strand diagram. (b) Rules for merges and splits.}
\label{fig:strand-diagram-circuit}
\end{figure}

We begin by describing the relationship between strand diagrams and piecewise-linear homeomorphisms for elements of Thompson's group~$F$.  Given a strand diagram for an element $f\in F$ and a real number $t\in[0,1]$, we can compute the image $f(t)$ using the procedure shown in Figure~\ref{fig:strand-diagram-circuit}(a).  Roughly speaking, the strand diagram acts like a computer circuit: whenever a number $t\in[0,1]$ is entered into the top, the signal winds its way through the circuit, emerging from the bottom as~$f(t)$.  The path that the signal takes is determined by its binary digits, which change as the signal passes through each node, as shown in Figure~\ref{fig:strand-diagram-circuit}(b).

\begin{figure}[b]
\centering
\hfill
\begin{minipage}[b]{80pt}
\centering
\includegraphics{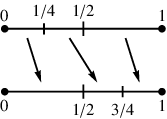}\\[-10pt]
\begin{align*}
.00\alpha &\mapsto .0\alpha \\
.01\alpha &\mapsto .10\alpha \\
.1\alpha &\mapsto .11\alpha
\end{align*}
\end{minipage}
\hfill\hfill
\includegraphics{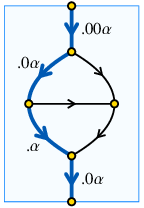}
\hfill
\includegraphics{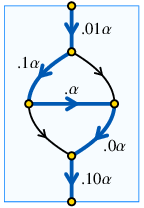}
\hfill
\includegraphics{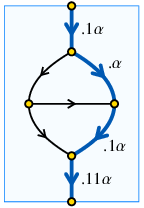}
\hfill
\caption{Three paths through a strand diagram.}
\label{fig:ThreePaths}
\end{figure}

Figure~\ref{fig:ThreePaths} shows the three different paths that a signal might take through
 the strand diagram for a certain element of~$F$.  Each of these paths corresponds to the linear mapping of an interval of the domain subdivision to the corresponding interval of the range subdivision.

\begin{figure}[b]
\centering
\includegraphics{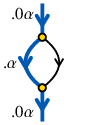}\;\;
\includegraphics{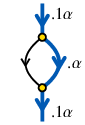}
\qquad\qquad\qquad
\includegraphics{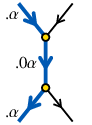}\;\;
\includegraphics{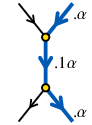}
\caption{Reductions do not change the underlying map.}
\label{fig:reductions-equivalent}
\end{figure}

Note that reductions do not change the action of the strand diagram on binary sequences, as shown in Figure~\ref{fig:reductions-equivalent}.  Thus equivalent strand diagrams really do represent the same piecewise-linear homeomorphism.

\begin{note}
More generally, we can interpret an $(m,n)$-strand diagram
as a \newword{Thompson-like homeomorphism}
$[0, m]\rightarrow[0, n]$, i.e. a piecewise-linear homeomorphism whose slopes
are powers of $2$, and whose breakpoints have dyadic rational coordinates.
Each number of the form $(k-1)+(.\alpha)$ corresponds to an input of .$\alpha$ entered into the $k$th source,
or an output of .$\alpha$ emerging from the $k$th sink.
%The set of Thompson-like
%homeomorphisms determined in this way is
%precisely the set of orientation-preserving homeomorphisms
%$[0, m]\rightarrow[0, n]$ whose slopes are powers of two, and whose breakpoints
%have dyadic rational coordinates.
\end{note}

%\begin{figure}
%\centering
%\subfloat[]{\includegraphics[height=2.3cm]{Three_Conjugate_Elements.eps}}
%\qquad\qquad\qquad\qquad
%\qquad \qquad
%\subfloat[]{\includegraphics[height=1.5cm]{Conjugate_Groupoid_Element.eps}}
%\qquad\qquad\qquad\qquad
%\qquad \qquad
%\subfloat[]{\includegraphics[height=2.3cm]{Groupoid_Example_Annulus.eps}}
%\caption{(a) Three conjugate elements (b) A minimal representative
%(c) The corresponding
%reduced annular strand diagram \blue{I want to remove this figure entirely.}}
%\label{fig:three-conjugate-minimal-annulus}
%\end{figure}

\subsection{Fixed points and replacement rules
\label{ssec:replacement-rules}
}

\begin{figure}[t]
\centering
\includegraphics{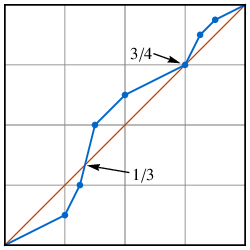}
\caption{An element of $F$.}
\label{fig:example-graph-F}
\end{figure}

Figure \ref{fig:example-graph-F} shows a typical element of Thompson's group~$F$.
This element has four fixed points at $0$, $1$, $1/3$, and~$3/4$, whose dynamics we now consider:
\begin{enumerate}
\item The fixed point at $0$ is attracting, since the slope is $1/2$.  Near the fixed point, the function acts on binary digits via the rule $f(.\alpha) = .0\alpha$.  This causes orbits to converge to~$0$:
\[
.\alpha \;\;\mapsto\;\; .0\alpha \;\;\mapsto\;\; .00\alpha \;\;\mapsto\;\; .000\alpha \;\;\mapsto\;
\; \cdots
\]
The fixed point at $1$ is similar, with $f(.\alpha) = .1\alpha$ in a neighborhood of~$1$.
\medskip

\item The fixed point at $1/3 = .\overline{10}\ldots$ is repelling, since the slope is $4$.  Near the fixed point, the function acts on binary digits via the rule $f(.10\alpha) = .\alpha$, which causes points to move away from~$1/3$:
\[
.101010\alpha \;\;\mapsto\;\; .1010\alpha \;\;\mapsto\;\; .10\alpha \;\;\mapsto\;\; .\alpha \;\;
\mapsto\;\; \cdots
\]

\item Since $3/4$ is a dyadic fraction, it has two binary expansions, namely $.10\overline{1}$ and $.11\overline{0}$.  Since the function has a breakpoint at $3/4$, the behavior is different on the two sides of the fixed point.  On the left side, the fixed point is attracting, with $f(.10\alpha) = .101\alpha$:
\[
.10\alpha \;\;\mapsto\;\; .101\alpha \;\;\mapsto\;\; .1011\alpha \;\;\mapsto\;\; .10111\alpha \;\;
\mapsto\;\; \cdots
\]
On the right side, the fixed point is repelling, with $f(.1100\alpha) = .110\alpha$:
\[
.110000\alpha \;\;\mapsto\;\; .11000\alpha \;\;\mapsto\;\; .1100\alpha \;\;\mapsto\;\; .110\alpha \;\;
\mapsto\;\; \cdots
\]
\end{enumerate}

%\subsection{\sout{Directed Loops and Fixed Points} \blue{Combine subsections.}
%\label{ssec:directed-loops-fixed-points}
%}

In general, an isolated fixed point for an element of $F$ is either repelling or attracting, though it may have two different behaviors if the fixed point is dyadic.  This behavior is evident from the action on binary digits: a repelling fixed point adds binary digits, while an attracting fixed point removes them.

Note that it is also possible for an element of~$F$ to have one or more closed intervals of fixed points.  We refer to a maximal interval of this form as a \textbf{fixed interval}.  Note that the endpoints of a fixed interval must always be dyadic.

It turns out that the behavior of the fixed points for an element $f\in F$ is invariant under conjugacy.  Indeed, it is possible to read this behavior directly from the annular strand diagram:

\begin{theorem}
\label{thm:directed-loops-fixed-point}
There is a correspondence between the fixed points of an element $f\in F$ and the closed loops in the corresponding reduced annular strand diagram.  In particular:
\begin{enumerate}
\item Every fixed interval corresponds to a free loop in the annular strand diagram.
\smallskip
\item Every isolated non-dyadic fixed point (as well as $0$ and $1$) corresponds to either a split loop or a merge loop.  In particular, a repelling fixed point with slope $2^n$ corresponds to a split loop with $n$ splits, and an attracting fixed point with slope $2^{-n}$ corresponds to a merge loop with $n$ merges.
\smallskip
\item Finally, every isolated dyadic fixed point (other than $0$ or $1$) corresponds to a pair of concentric loops, which may be split loops, merge loops, or one of each, depending on the behavior on the two sides of the fixed point.
\end{enumerate}
In the latter two cases, the pattern of outward and inward connections around the loop determines
the tail of the binary expansion of the fixed point.
Specifically, each outward connection corresponds to a $1$, and each inward connection corresponds to a $0$.
\end{theorem}

\begin{proof}
By Proposition~\ref{prop:ReductionConjugate}, there exists an $(n,n)$-strand diagram $\mathfrak{g}$ whose closure is the reduced annular strand diagram for~$f$.  This diagram $\mathfrak{g}$ corresponds to a Thompson-like homeomorphism $g\colon[0,n]\to[0,n]$ which is conjugate to~$f$.  The fixed points of $g$ are in one-to-one correspondence with the fixed points of~$f$, and indeed it suffices to prove the statement of the theorem for the homeomorphism~$g$.

To begin, suppose that the reduced annular strand diagram for $g$ contains a merge loop, e.g.~the loop shown in Figure~\ref{fig:example-merge-loop}(a).  In the strand diagram for $g$, this loop corresponds to a directed path from the $k$th source to the $k$th sink, all of whose vertices are merges.  In particular, if we feed a binary number $.\alpha$ into the $k$th source, it will emerge from the $k$th sink with some finite prefix added.  For example, Figure~\ref{fig:example-merge-loop}(b) shows how an input of $.\alpha$ into the $k$th source will become output of $.1101\alpha$ from the $k$th sink.  In terms of functions, this means that
\[
g\bigl((k-1)+.\alpha)\bigr) = (k-1)+(.1101\alpha)
\]
for all $.\alpha\in [0,1]$.  It follows that $(k-1) + .\overline{1101}$ is an attracting fixed point for~$g$.

In this way, each merge loop corresponds to an attracting fixed point of~$g$, and each split loop corresponds to a repelling fixed point for $g$.  A free loop corresponds to an entire unit interval $[k,k+1]$ of fixed points.  Note that every fixed point arises in this fashion.  In particular, every fixed point must lie in some interval $[k,k+1]$, and such an interval contains a fixed point if and only if there is a directed path in the strand diagram from the $k$th source to the $k$th sink.  Also, note that each interval $[k,k+1]$ contains at most one fixed point, and such a fixed point is dyadic if and only if it is an endpoint of the interval.

Finally, observe that the conjugacy between $f$ and $g$ preserves which fixed points are dyadic, and also preserves the tail of the binary expansion of each fixed point.

\begin{figure}
\centering
\subfloat[]{\includegraphics{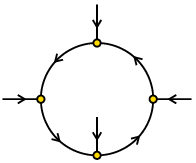}}
%\qquad\qquad\qquad\qquad
\qquad \qquad
\qquad
\subfloat[]{\includegraphics{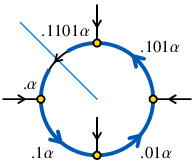}}
\caption{(a) An example of a merge loop (b)
Traveling around the merge loop}
\label{fig:example-merge-loop}
\end{figure}
\end{proof}

Note that the outermost loop of an annular strand diagram for an element of $F$ corresponds to the fixed point
$0=.0000\cdots$, while the innermost loop corresponds to the fixed point $1=.1111\cdots$.
Within each connected component, the outermost and innermost loops correspond to dyadic fixed points, while the interior loops correspond to non-dyadic fixed points.

\begin{corollary}
Let\/ $\mathfrak{f}'$ be the reduced annular strand diagram for an element $f\in F$.
Then every component of\/ $\mathfrak{f}'$ corresponds to exactly one of the following:
\begin{enumerate}
\item A maximal closed interval of fixed points of $f$ (for a free loop), or\smallskip
\item A maximal interval with no dyadic fixed points of $f$ in its interior.
\end{enumerate}
\end{corollary}

\begin{figure}[b]
\centering
\hfill
\includegraphics{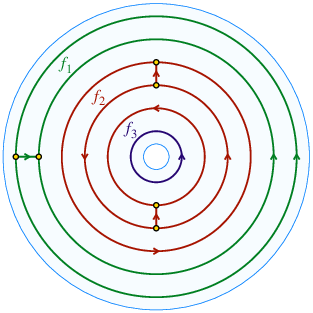}
\hfill\hfill
\includegraphics{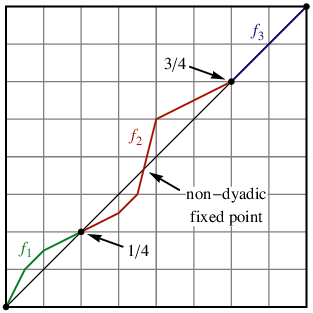}
\hfill
\caption{An element of $F$ with three components, and the corresponding annular strand diagram.}
\label{fig:components}
\end{figure}

If $f\in F$, a \newword{cut point} of $f$ is either an isolated dyadic fixed point of~$f$,
or an endpoint of a maximal interval of fixed points.
If $0 = \alpha_0 < \alpha_1 < \cdots < \alpha_n = 1$ are the cut points of $f$, then the restrictions
$f_i\colon [\alpha_{i-1}, \alpha_i] \rightarrow [\alpha_{i-1}, \alpha_i]$ are called the \newword{components}
of~$f$ (see Figure~\ref{fig:components}).
Each component of $f$ corresponds to one connected component of the reduced annular strand diagram
(figure \ref{fig:components}).
If $\alpha < \beta$ are any dyadic rationals, it is well known (see~\cite{cfp})
that there exists a Thompson-like homeomorphism
$\varphi\colon [\alpha, \beta]\rightarrow [0,1]$.
It follows that any Thompson-like homeomorphism of $[\alpha, \beta]$ can be conjugated
by $\varphi$ to give an element of $F$. The following are straightforward.

\begin{proposition}
Let $f\in F$ have components $f_i\colon [\alpha_{i-1}, \alpha_i] \rightarrow [\alpha_{i-1}, \alpha_i]$,
and let\/ $\mathfrak{f}'$ be the reduced annular strand diagram for $f$.  Then for each $i$, the component of $S$ corresponding to $f_i$
is the reduced annular strand diagram for any element of $F$ conjugate to $f_i$.
\end{proposition}

%\begin{proof}
%Suppose $f$ has $n + 1$ cut points $0 = \alpha_0 < \alpha_1 < \cdots < \alpha_n = 1$.
%Then we can conjugate $f$ to an
%element of Thompson's groupoid whose cut points are at $0, 1, 2, \ldots, n$.  The %resulting $(n, n)$-strand diagram
%has $n$ connected components which, when reduced, yield the $n$ components of $S$.
%\qed
%\end{proof}

\begin{corollary}
Let $f, g\in F$ have components $f_1, \ldots, f_n$ and $g_1, \ldots, g_n$.  Then $f$ is conjugate to $g$
in $F$ if and only if each $f_i$ is conjugate to $g_i$ through some Thompson-like homeomorphism.
\end{corollary}

\begin{remark}
Although the work outlined in Sections \ref{sec:representation-elements} and
\ref{ssec:replacement-rules} describes a relation between
the dynamics of elements of $F$ and annular strand diagrams, one can describe
a similar connection between elements in $T$ (respectively, $V$) and
toral strand diagrams  (respectively, abstract closed strand diagrams).
Bleak \emph{et al.}
give a description of centralizers for Thompson's group $V$
based on a study of the orbits of an element~\cite{matucci8}.
The description in~\cite{matucci8} recovers abstract closed strand diagrams and
shows how they encode information about the dynamics of all the elements
of a conjugacy class, giving a generalization of the results obtained above for $F$.
\end{remark}

%%%%%%%%%%%%%%%%%%%%%%%%%%%%%%%%%%%%%
%%%%%%%%%%%%%%%%%%%%%%%%%%%%%%%%%%%%%
%%%%%%%%%%%%%%%%%%%%%%%%%%%%%%%%%%%%%
%%%%%%%%%%%%%%%%%%%%%%%%%%%%%%%%%%%%%
%%%%%%%%%%%%%%%%%%%%%%%%%%%%%%%%%%%%%
%%%%%%%%%%%%%%%%%%%%%%%%%%%%%%%%%%%%%
%%%%%%%%%%%%%%%%%%%%%%%%%%%%%%%%%%%%%
%%%%%%%%%%%%%%%%%%%%%%%%%%%%%%%%%%%%%
%%%%%%%%%%%%%%%%%%%%%%%%%%%%%%%%%%%%%
%%%%%%%%%%%%%%%%%%%%%%%%%%%%%%%%%%%%%
%%%%%%%%%%%%%%%%%%%%%%%%%%%%%%%%%%%%%
%%%%%%%%%%%%%%%%%%%%%%%%%%%%%%%%%%%%%
%%%%%%%%%%%%%%%%%%%%%%%%%%%%%%%%%%%%%
%%%%%%%%%%%%%%%%%%%%%%%%%%%%%%%%%%%%%
%%%%%%%%%%%%%%%%%%%%%%%%%%%%%%%%%%%%%
%%%%%%%%%%%%%%%%%%%%%%%%%%%%%%%%%%%%%
%%%%%%%%%%%%%%%%%%%%%%%%%%%%%%%%%%%%%
%%%%%%%%%%%%%%%%%%%%%%%%%%%%%%%%%%%%%
%%%%%%%%%%%%%%%%%%%%%%%%%%%%%%%%%%%%%
%%%%%%%%%%%%%%%%%%%%%%%%%%%%%%%%%%%%%
%%%%%%%%%%%%%%%%%%%%%%%%%%%%%%%%%%%%%
%%%%%%%%%%%%%%%%%%%%%%%%%%%%%%%%%%%%%
%%%%%%%%%%%%%%%%%%%%%%%%%%%%%%%%%%%%%
%%%%%%%%%%%%%%%%%%%%%%%%%%%%%%%%%%%%%
%%%%%%%%%%%%%%%%%%%%%%%%%%%%%%%%%%%%%
%%%%%%%%%%%%%%%%%%%%%%%%%%%%%%%%%%%%%
%%%%%%%%%%%%%%%%%%%%%%%%%%%%%%%%%%%%%
%%%%%%%%%%%%%%%%%%%%%%%%%%%%%%%%%%%%%

\subsection{Mather Invariants \label{sec:mather-invariant}}

Conjugacy in $F$ was first investigated by Brin and Squier \cite{brin2}, who successfully
found an invariant for conjugacy in the full group of piecewise-linear
homeomorphisms of the interval with finitely many breakpoints.
Their invariant was based on some ideas of
Mather \cite{mather} for determining whether two given diffeomorphisms
of the unit interval are conjugate. In this section we show that the annular strand diagrams described in Section~\ref{sec:unified-pow} can be used to define a Mather-type invariant
for elements of $F$ for functions without fixed points (except for $0$ and $1$).
The invariant that we describe is a suitable adaptation of Brin and Squier
which works for $F$.
independently from our work, Gill and Short \cite{ghisho1}
also successfully found a way to generalize Brin and Squier's version of the Mather
invariant and thus providing a description of conjugacy in $F$ which is similar
to the one presented in the current section (although the proofs are somewhat
different).

Recall that a map $f \in F$ is called a \textbf{one-bump function}
is an element such that $f(x) > x$ for all $x \in (0, 1)$.
Consider a one-bump function $f\in F$, with slope $2^{m_0}$ at $0$ and
slope $2^{m_1}$ at $1$.  In a neighborhood of zero,
$f$ acts as multiplication by $2^{m_0}$;
in particular, for any sufficiently small $t > 0$, the interval $[t, 2^{m_0}t]$
is a fundamental domain for the action of $f$.
%(see figure \ref{fig:mather-circle-identification}).
%\begin{figure}[0.5\textwidth]
%\includegraphics{Mather_Circle_Identification.eps}
%\centering
%\caption{Action of $f$ in a neighborhood of $0$}
%\label{fig:mather-circle-identification}
%\end{figure}
If we make the identification $t \sim 2^{m_0}t$ in the interval $(0, \epsilon)$, we obtain a circle $C_0$, with partial covering map
$p_0 \colon (0, \epsilon)\rightarrow C_0$.
%Note that the restriction of $f$ is a deck transformation of this cover:
%$$
%\begindc{\commdiag}[12]
%\obj(0,40){$(0,\epsilon)$}
%\obj(60,40){$(0, \epsilon)$}
%\obj(30,0){$C_0$}
%\mor{$(0,\epsilon)$}{$(0, \epsilon)$}{$f$}
%\mor{$(0,\epsilon)$}{$C_0$}{$p_0$}[\atright, \solidarrow]
%\mor{$(0, \epsilon)$}{$C_0$}{$p_0$}[\atleft, \solidarrow]
%\enddc
%$$
Similarly, if we identify $(1 - t) \sim (1 - 2^{m_1} t)$ on the interval $(1 - \delta, 1)$, we obtain a circle $C_1$, with
partial covering map $p_1\colon (1 - \delta, 1)\rightarrow C_1$.

If $N$ is sufficiently large, then $f^N$ will take some lift of $C_0$ to $(0, \epsilon)$
and map it to the interval $(1 - \delta, 1)$.
This induces a map $f^\infty\colon C_0\rightarrow C_1$, making the following diagram commute:
$$
\begindc{\commdiag}[1]
\obj(-50,550){$(0, \epsilon)$}
\obj(700,550){$(1-\delta, 1)$}
\obj(0,0){$C_0$}
\obj(550,0){$C_1$}
\mor(0,550)(550,550){$f^N$}[\atleft, \solidarrow]
\mor(0,550)(0,0){$p_0$}[\atright, \solidarrow]
\mor(550,550)(550,0){$p_1$}[\atleft, \solidarrow]
\mor(0,0)(550,0){$f^{\infty}$}[\atleft, \solidarrow]
\enddc
$$
%\definition{The map $f^\infty$ defined above is the Mather invariant for $f$.}
We note that $f^\infty$ does not depend on the specific value of $N$ chosen. Any map $f^s$,
for $s \ge N$, induces the same map $f^\infty$. This is because $f$ acts as the identity on $C_1$ by construction
and $f^s$ can be written as $f^{s-N}(f^N(t))$, with $f^N(t) \in (1-\delta,1)$.
%If $k > 0$, then the map $t \mapsto kt$ on $(0, \epsilon)$ induces
%a ``rotation'' $\rot_k$ of $C_0$.
%In particular, if we use the
%coordinate $\theta = \log t$ on $C_0$, then
%\begin{equation*}
%\rot_k (\theta) \,=\, \theta + \log k
%\end{equation*}
%so rot$_k$ is an actual rotation.	

\begin{remark}
The map $f^\infty$ defined above is the \textbf{Mather invariant} for $f$ as defined
in Brin and Squier \cite{brin2} for functions in
the full group of piecewise-linear homeomorphisms of the interval with
finitely many breakpoints. Since $f$ is an element of $F$,
we will now rescale the two circles $C_0$ and $C_1$
so that their length is equal to the exponents of the slopes at $0$ and $1$.
\end{remark}

\begin{definition}
The \textbf{piecewise-linear logarithm} $\PLog\colon (0, \infty)\rightarrow (-\infty, \infty)$
is the piecewise-linear function that maps the interval $\left[ 2^k, 2^{k+1} \right]$ linearly onto
$[k, k+1]$ for every $k \in \Zb$, as shown below:
\begin{figure}[h]
\centering
\includegraphics[height=2cm]{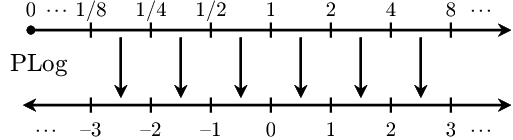}
%\caption{The $\PLog$ map}
%\label{fig:PLog-map}
\end{figure}
\end{definition}
Suppose that $f\in F$ is a one-bump function with slope $2^m$ at
$0$ and slope $2^{-n}$ at $1$,
and let $f^\infty\colon C_0 \rightarrow C_1$ be the corresponding Mather invariant.
In a neighborhood of $0$, the function $f$ acts as multiplication by $2^m$.  In particular,
$\PLog f(t) = m + \PLog t$ for all $t \in (0,\epsilon)$, so we can identify $C_0$ with the
circle $\Rb / m \Zb$.  The following picture shows the case $m = 3$:
\begin{figure}[h]
\centering
\includegraphics[height=2.5cm]{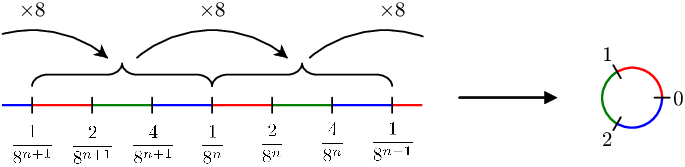}
%\caption{Construction of the circle $C_0$}
\label{fig:dyadic-circle-identification}
\end{figure}

In a similar way, we can use the function $t \mapsto -\PLog(1- t)$ to identify $C_1$ with the circle $\Rb / n \Zb$.
This lets us regard the Mather invariant for $f$ as a function $f^\infty\colon \Rb /m\Zb \rightarrow \Rb /n \Zb$.
Because $f^N$ and $\PLog$ are piecewise-linear, the Mather invariant $f^\infty$ is a piecewise-linear function.
Moreover, $f^\infty$ is Thompson-like: all the slopes are powers of $2$, and the breakpoints are dyadic rational
numbers of $\Rb/m \Zb = [0,m]/\{0,m\}$.

Now, if $k \in \Zb$, then the map $t \mapsto 2^kt$ on $(0, \epsilon)$ induces precisely an
integer rotation of $\Rb / m \Zb$:
\begin{equation*}
\rot_k (\theta) = \theta + k \mod m
\end{equation*}
We are now ready to state a criterion for conjugacy via Mather invariants:

\begin{theorem}
\label{MainMatherTheorem}
Let $f, g \in F$ be one-bump functions with $f'(0) = g'(0) = 2^m$ and $f'(1) = g'(1) = 2^{-n}$,
and let $f^\infty, g^\infty\colon \Rb / m \Zb \rightarrow \Rb / n \Zb$ be the corresponding Mather invariants.
Then $f$ and $g$ are conjugate if and only if $f^\infty$ and $g^\infty$ differ by integer rotations of the domain and range circles:
\begin{equation*}
\begindc{\commdiag}[17]
\obj(0,30) {$\Rb / m\Zb$}
\obj(40,30){$\Rb / n\Zb$}
\obj(0,0){$\Rb / m   \Zb$}
\obj(40,0){$\Rb / n   \Zb$}
\mor{$\Rb / m\Zb$}{$\Rb / n\Zb$}{$f^\infty$}[\atleft, \solidarrow]
\mor{$\Rb / m   \Zb$}{$\Rb / n   \Zb$}{$g^\infty$}[\atleft, \solidarrow]
\mor{$\Rb / m\Zb$}{$\Rb / m   \Zb$}{$\rot_k$}[\atright, \solidarrow]
\mor{$\Rb / n\Zb$}{$\Rb / n   \Zb$}{$\rot_\ell$}[\atleft, \solidarrow]
\enddc
\end{equation*}
\end{theorem}

We only sketch the proof of Theorem \ref{MainMatherTheorem}
as some steps are straightforward,
while others can be obtained using small variations of some of our previous arguments.

The forward direction follows observing that, if $f = h^{-1} g h$ for some
$h \in F$, then the following diagram commutes, where $k = \log_2 h'(0)$ and
$\ell = \log_2 h'(1)$:
\begin{equation*}
\begindc{\commdiag}[20]
\obj(0,0){$\Rb / m \Zb \; \; \; \; \;$}
\obj(50,0){$\; \; \; \; \; \Rb / n \Zb$
}
\obj(20,20){$\Rb / m \Zb \; \; \; \; \; $}
\obj(70,20){$ \; \; \; \; \; \Rb / n \Zb$}
\obj(0,50){$(0,\epsilon)$}
\obj(50,50){$(1-\delta,1)$}
\obj(20,70){$(0, \epsilon)$}
\obj(70,70){$(1-\delta, 1)$}
\mor(0,0)(50,0){$f^\infty$}[\atright, \solidarrow]
\mor(0,0)(20,20){$\rot_k$}[\atright, \solidarrow]
\mor(50,0)(70,20){$\rot_\ell$}[\atright, \solidarrow]
\mor(20,20)(50,20){$\quad g^\infty$}[\atleft, \solidline]
\mor(50,20)(70,20){}[\atleft, \solidarrow]

\mor(0,50)(0,0){$p_0$}[\atright, \solidarrow]
\mor(20,70)(20,50){}[\atleft, \solidline]
\mor(20,50)(20,20){$p_0$}[\atleft, \solidarrow]
\mor(50,50)(50,0){$\begin{matrix} p_1 \\ \, \end{matrix}$}[\atleft, \solidarrow]
\mor(70,70)(70,20){$p_1$}[\atleft, \solidarrow]

\mor{$(0,\epsilon)$}{$(1-\delta,1)$}{$\qquad f^N$}[\atleft, \solidarrow]
\mor{$(0,\epsilon)$}{$(0, \epsilon)$}{$h$}[\atleft, \solidarrow]
\mor{$(1-\delta,1)$}{$(1-\delta, 1)$}{$h$}[\atleft, \solidarrow]
\mor{$(0, \epsilon)$}{$(1-\delta, 1)$}{$g^N$}[\atleft, \solidarrow]
\enddc
\end{equation*}

For the converse, we must show that any
two one-bump functions whose Mather invariants
differ by integer rotation are conjugate in $F$.  To prove this, we describe an explicit correspondence between
Mather invariants and reduced annular strand diagrams.

\begin{figure}[t]
\centering
\subfloat[]{\includegraphics{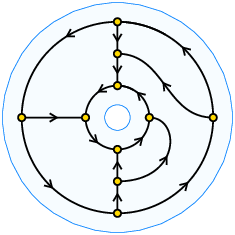}}
\qquad\qquad\qquad\qquad
\subfloat[]{\includegraphics{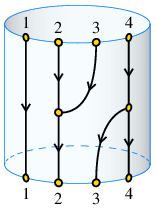}}
\caption{(a) Annular strand diagram for a one-bump function (b) The corresponding cylindrical strand diagram}
\label{fig:annulus-to-cylinder}
\end{figure}

If $f\in F$ is a one-bump function, then the only fixed points of $f$ are at $0$ and $1$.  Therefore, the reduced annular
strand diagram for $f$ has only two directed cycles (see Figure~\ref{fig:annulus-to-cylinder}a).
Since $f'(0) > 1$, the outer cycle (corresponding to $0$) must be a split loop, and the inner cycle (corresponding to $1$)
must be a merge loop.  If we remove these two cycles, we get an  $(m,n)$-cylindrical
strand diagram (see Figure~\ref{fig:annulus-to-cylinder}b). Such a diagram can be used to describe a Thompson-like map between two circles.
The following result only requires one to observe how to construct and glue
forests to obtain a cylindrical strand diagram and conversely how to cut
such diagram to obtain a pair of forests. We omit its proof.

\begin{proposition}There is a one-to-one correspondence between
\begin{enumerate}
\item Reduced cylindrical $(m,n)$-strand diagrams, and
\item Thompson-like functions $\Rb / m \Zb \rightarrow \Rb / n \Zb$, with two functions considered equivalent if they differ
by integer rotation of the domain and range circles.
\end{enumerate}
\end{proposition}

To complete the proof of Theorem \ref{MainMatherTheorem}
we need the following result which can be obtained by imitating the idea of
Theorem \ref{thm:directed-loops-fixed-point}
to study the image of a binary number under many iterations.

\begin{proposition}
Let $\mathcal{A}$ be the reduced annular strand diagram for a one-bump function $f \in F$, and let
$\mathcal{C}$ be the cylindrical $(m,n)$-strand diagram obtained by removing the merge and split loops from $\mathcal{A}$.
Then $\mathcal{C}$ is the cylindrical strand diagram for the Mather invariant $f^\infty\colon \Rb /m\Zb \rightarrow \Rb /n\Zb$.
\end{proposition}

\begin{remark} Brin and Squier used the Mather invariant as a means to describe
a conjugacy invariant in the full group of piecewise-linear homeomorphisms
of the unit interval with finitely many breakpoints. With
Theorem~\ref{MainMatherTheorem}, the authors found a way to adapt
Brin and Squier's invariant to the case of Thompson's group $F$.
After the authors
found the description of Theorem~\ref{MainMatherTheorem}
and proved it using strand diagrams, the second author analyzed this description
using the techniques appearing in the joint paper of Kassabov and the second author
\cite{matucci5}: this resulted in a generalization \cite{matucci3}
which works for a larger class of
groups of piecewise-linear homeomorphisms (essentially the Thompson-Stein groups)
and bridging the gap between the descriptions in \cite{brin2} and Theorem~\ref{MainMatherTheorem}.
\end{remark}

\section*{Acknowledgments}
The authors would like to thank Ken Brown and Martin Kassabov for a very careful reading of this paper
and many helpful remarks.
The authors would also like to thank Stephen Pride for providing fundamental references and
Collin Bleak, Matt Brin, J\"org Lehnert and Mark Sapir for helpful conversations.
Finally, the authors would like to
thank the referee for many helpful comments and suggestions.

\bibliographystyle{plain}
\bibliography{go4}

\end{document}